\documentclass[12pt]{amsart}

\usepackage{graphics, amsmath, amsfonts, amssymb, amscd, mathrsfs}

\newtheorem{theorem}{Theorem}

\theoremstyle{definition}

\newcommand{\C}{\mathbb{C}}
\newcommand{\codim}{\operatorname{codim}}

\title{Smooth toric varieties are Oka}

\author{Finnur L\'arusson}
\address{School of Mathematical Sciences, University of Adelaide, Adelaide SA 5005, Australia.} 
\email{finnur.larusson@adelaide.edu.au}



\date{1 July 2011.  Latest minor changes 13 October 2013.}

\begin{document}

\begin{abstract}
We show that every smooth toric variety over the field of complex numbers is an Oka manifold.
\end{abstract}

\maketitle

\noindent
The class of Oka manifolds has emerged from the modern theory of the Oka principle, initiated in 1989 in a seminal paper of Gromov \cite{Gromov}.  They were first formally defined by Forstneri\v c in 2009 in the wake of his result that some dozen possible definitions are all equivalent \cite{Forstneric2}.  Often, the Oka property is established by verifying the sufficient geometric condition of ellipticity, introduced by Gromov.  The Oka property can be seen as an answer to the question: what should it mean for a complex manifold to be \lq\lq anti-hyperbolic\rq\rq?  For more background, see the survey \cite{Forstneric-Larusson}.

The first examples of Oka manifolds are complex Lie groups and their homogeneous spaces.  Further examples are listed in \cite{Forstneric-Larusson}, Section 3.  In this note, we show that all smooth toric varieties over the field of complex numbers are Oka.  For the theory of toric varieties we refer to the monograph \cite{Cox-Little-Schenck}, primarily Section 5.1.

The facts from Oka theory that we need are the following.

\begin{theorem}[\cite{Gromov}, 0.5.B]
\label{examples}
A complex Lie group is elliptic and thus Oka.  The complement in $\C^n$ of an algebraic subvariety of codimension at least $2$ is elliptic and thus Oka.
\end{theorem}

\begin{theorem}[\cite{Forstneric1}, Theorem 1.4; \cite{Forstneric2}, Corollary 1.3]
\label{up-and-down}
If $E$ and $B$ are complex manifolds and $E\to B$ is a holomorphic fibre bundle
whose fibres are Oka manifolds, then $B$ is an Oka manifold if and only if $E$ is an Oka manifold. 
\end{theorem}

Now let $X$ be a smooth toric variety over $\C$.  If $X$ has a torus factor, say $X$ is isomorphic to $Y\times(\C^*)^k$, $k\geq 1$, where $Y$ is another smooth toric variety, then, by Theorems \ref{examples} and \ref{up-and-down}, $X$ is Oka if and only if $Y$ is Oka.  

Hence we may assume that $X$ has no torus factor, so the quotient construction described in \cite{Cox-Little-Schenck}, Section 5.1, applies; see in particular Theorem 5.1.11.  We can write $X$ as a geometric quotient
\[ X = (\C^m\setminus Z)/G, \]
where $Z$ is a union of coordinate subspaces of $\C^m$, and the group $G$ is a subgroup of $(\C^*)^m$ acting on $\C^m\setminus Z$ by diagonal matrices.  In fact, $G$ is isomorphic to the product of a torus and a finite abelian group (\cite{Cox-Little-Schenck}, Lemma 5.1.1), so $G$ is reductive.  Also, $G$ is Oka (Theorem \ref{examples}).  Furthermore, $\codim Z\geq 2$ (\cite{Cox-Little-Schenck}, Exercise 5.1.13), so $\C^m\setminus Z$ is Oka (Theorem \ref{examples}).

Since $X$ is smooth, $G$ acts freely on $\C^m\setminus Z$ (\cite{Cox-Little-Schenck}, Exercise 5.1.11).  We would like to conclude that the projection $\C^m\setminus Z\to X$ is a holomorphic fibre bundle.  By Theorem \ref{up-and-down}, this would imply that $X$ is Oka.

Note that $Z$, being a union of coordinate subspaces, is the intersection of unions of coordinate hyperplanes.  Thus $\C^m\setminus Z$ is the union of Zariski-open sets of the form $U=\C^m\setminus (H_1\cup\cdots\cup H_k)$, where $H_1,\dots,H_k$ are coordinate hyperplanes.  Each $U$ is affine algebraic or, from the holomorphic point of view, Stein, as well as $G$-invariant.  By slice theory for actions of reductive groups, the quotient map $U\to U/G$ is a holomorphic fibre bundle (\cite{Snow}, Corollary 5.5), or, from the algebraic point of view, a locally trivial fibration in the \'etale sense (\cite{Luna}, Corollaire 5).  It follows that $\C^m\setminus Z\to X$ is a holomorphic fibre bundle.

Thus we have proved the following result.

\begin{theorem}
\label{our-theorem}
Every smooth toric variety over $\C$ is Oka.
\end{theorem}

We conclude with three remarks.  

(a)  It is not clear from the above whether a smooth toric variety is elliptic, even though $\C^m\setminus Z$ is.

(b)  By definition, every toric variety is rational, that is, birationally equivalent to complex projective space, which is Oka.  Thus, if we knew that the Oka property was birationally invariant, our result would be immediate.  At present, it is not even known how the Oka property behaves with respect to blowing up a point, except in some special cases.

(c)  It follows from \cite{AKZ}, Theorem 2.1, which is much more difficult to prove than our result, that the smooth locus of an \textit{affine} toric variety is elliptic and thus Oka.  (Smooth affine toric varieties are of the form $\C^n \times (\C^*)^m$, so they are obviously elliptic.)

\end{document}